\documentclass[12pt]{amsart}
\usepackage{amssymb}
\usepackage{amsthm}

\usepackage{amscd}
\usepackage{extarrows}
\usepackage{graphicx}
\usepackage{fancyhdr}
\usepackage{epsfig}
\usepackage{amsfonts}
\usepackage{mathrsfs}
\usepackage{picinpar} \usepackage{amsmath} \usepackage[all]{xy}
\usepackage[colorlinks=true]{hyperref}
\hypersetup{urlcolor=blue, citecolor=red}
\theoremstyle{plain}
\newtheorem{theorem}{Theorem}[section]

\newtheorem{lemma}[theorem]{Lemma}
\newtheorem{proposition}[theorem]{Proposition}
\theoremstyle{definition}
\newtheorem{definition}[theorem]{Definition}

\begin{document}
\title{Gorenstein homological dimensions of modules over triangular matrix rings}
\author{Rongmin Zhu, Zhongkui Liu, and Zhanping Wang}
\footnote[0]{2014.12.30}
\date{} \maketitle

\hspace{6.3cm}\noindent{\footnotesize {\bf Abstract}
\vspace{0.2cm}

\hspace{-0.45cm}Let $A$ and $B$ be rings, $U$ a $(B, A)$-bimodule and $T=\left(\begin{smallmatrix} A & 0 \\  U & B \\\end{smallmatrix}\right)$ be the triangular matrix ring. In this paper, we characterize the Gorenstein
homological dimensions of modules over $T$, and discuss when a left $T$-module is strongly Gorenstein projective or strongly Gorenstein injective module.
\vspace{0.2cm}

\noindent{\footnotesize {{\it{Keywords}}:} Triangular matrix ring; Gorenstein regular ring; Gorenstein homological dimension

\section{Introduction and Preliminaries}
Triangular matrix rings have been studied by many authors (e.g. see \hyperref[14-16]{[14--16]} and their references). Such rings play an important role in the study of the representation theory of Artin rings and algebras. The modules (left or right) over such rings can be
described in a very concrete fashion and we have nice descriptions of some important
classes of modules over such rings. Krylov and Tuganbaev [18] presented general properties of matrix rings and injective, projective, flat, and hereditary modules over such rings.
Let $A$ and $B$ be rings and $U$ a $(B,A)$-bimodule. We denote by $T$ the triangular matrix ring $\left(\begin{smallmatrix}  A & 0 \\  U & B \\\end{smallmatrix}\right)$. Using the description of $T$-modules, Asadollahi and Salarian [1] studied the vanishing of the extension functor ¡®Ext¡¯
over $T$ and explicitly described the structure of (right) $T$-modules of finite
projective (resp. injective) dimension.  Enochs and Torrecillas [4] described flat covers and cotorsion envelopes of modules
over $T$.

In the 1990s Enochs, Jenda, and Torrecillas introduced the Gorenstein projective, injective, and flat modules [6, 8] and then developed Gorenstein homological algebra [7]. Zhang [19] introduced in 2013 the compatible bimodules and explicitly described the Gorenstein projective modules over triangular matrix Artin algebra. In 2014, Enochs and other authors in [5] introduced Gorenstein regular ring and characterized when a left $T$-module is Gorenstein projective or Gorenstein injective over such rings.

This paper is devoted to study Gorenstein homological dimensions over triangular matrix rings and organized as follows. In Section 2, we focus on discussing the structure
 of left $T$-modules of finite Gorenstein projective (resp. injective) dimensions. Let $n$ be a non-negative integer, $B$ a left Gorenstein regular ring, $U_{A}$ have finite flat dimension and $_{B}U$ have finite projective dimension. Let
$M=\binom{M_{1}}{M_{2}}_{\varphi^{M}}$  be a left $T$-module. Then using the structure of left $T$-modules, we show that $\mathrm{Gpd}(M)\leq n$ if and only if $\mathrm{Gpd}(M_{1})\leq n$,
$\mathrm{Gpd}(\frac{M_{2}}{\mathrm{Im}\varphi^{M}})\leq n$ and if
$\binom{K_{1}}{K_{2}}_{\varphi^{K}}$ is a n-th syzygy of $M$, then
$\varphi^{K}$ is injective, where $\mathrm{Gpd}(M)$ denotes the Gorenstein projective dimension of $M$ (Theorem 2.5). Similar dual result holds for Gorenstein injective dimension (Theorem 2.6).

Motivated by the characterization of Gorenstein projective or Gorenstein injective left $T$-module, in Section 3, we characterize when a left $T$-module is Gorenstein flat. We prove in Theorem 3.9 that if $U_{A}$ and $_{B}U$ are finitely generated and have finite projective dimension, $T$ is a Gorenstein ring, then $F=\binom{F_{1}}{F_{2}}_{\varphi^{F}}\in\mathcal{GF}(T)$ if and only if $F_{1}\in\mathcal{GF}(A)$, $\frac{F_{2}}{\mathrm{Im}\varphi^{F}}\in\mathcal{GF}(B)$ and $\varphi^{F}$ is injective, where $\mathcal{GF}$ denotes the class of all Gorenstein flat modules.

There is also an analogoue for free module, namely, the strongly Gorenstein projective
module [3]. As observed by Bennis and Mahdou, a module is Gorenstein projective
if and only if it is a direct summand of a strongly Gorenstein projective module ([3,
Theorem 2.7]). Gao and Zhang [11] gave a concrete construction of strongly Gorenstein projective
modules, via the existed construction of upper triangular matrix Artin algebras. Finally, we give in Section 4 the characterization of strongly Gorenstein projective (resp. injective) modules and dimensions, which extend the results in [11]. We show in Theorem 4.2 that: if $U_{A}$ has finite flat dimension and $_{B}U$ has finite projective dimension, $B$ is left Gorenstein regular, then $M=\binom{M_{1}}{M_{2}}_{\varphi^{M}}$ is strongly Gorenstein projective if and only if $M_{1}$ and $\frac{M_{2}}{\mathrm{Im}\varphi^{M}}$ are strongly Gorenstein projective and the $\varphi^{M}$ is injective.

Throughout this paper, all rings are associative rings with identity, and all modules are unitary. As usual, $\mathrm{pd}(M),~\mathrm{id}(M)$ and $\mathrm{fd}(M)$ denote the projective, injective and flat dimensions of a left $R$-module $M$, respectively.

Let $R$ be a ring. A left $R$-module $M$ is called $Gorenstein~projective$ if there exists an exact sequence $$ \cdots \longrightarrow P^{-2} \longrightarrow  P^{-1} \longrightarrow  P^{0} \longrightarrow  P^{1} \longrightarrow  \cdots$$
of projective left $R$-modules with $M\cong\mathrm{Ker}(P^{0}\rightarrow P^{1})$ such that $\mathrm {Hom}_{R}(-,Q)$ leaves the sequence exact for any projective left $R$-module $Q$. The Gorenstein injective modules are defined dually. A left $R$-module $M$ is said to be $Gorenstein~flat$ if there exists an exact sequence $$ \cdots \longrightarrow F^{-2} \longrightarrow  F^{-1} \longrightarrow  F^{0} \longrightarrow  F^{1} \longrightarrow  \cdots$$
of flat left $R$-modules with $M\cong\mathrm{Ker}(F^{0}\rightarrow F^{1})$ such that $I\otimes_{R}-$ leaves the sequence exact for any injective right $R$-module $I$. We denote by  $\mathcal{GP}(R), ~\mathcal{GI}(R),$ and $~\mathcal{GF}(R)$ the class of all Gorenstein projective, injective and flat $R$-modules, respectively. For any $R$-module $M$, $\mathrm{Gpd}(M),~\mathrm{Gid}(M),$ and $\mathrm{Gfd}(M)$ denote the Gorenstein projective, injective and flat dimension of $M$, respectively, and $\mathrm{glGpd}(R)$ and $\mathrm{glGid}(R)$ denote the global Gorenstein projective and injective dimensions of $R$, respectively.

A complex $\mathbf{C}$ of modules is a sequence $$\cdots \longrightarrow C^{n-1}\stackrel{d^{n-1}}\longrightarrow C^{n}\stackrel{d^{n}}\longrightarrow C^{n+1}\stackrel{}\longrightarrow\cdots$$ of $R$-modules and $R$-homomorphisms such that $d^{n}d^{n-1}=0$ for all $n\in \mathbb{Z}$. A complex $\mathbf{C}$ is exact if for each $n$, $\mathrm{Ker}d^{n}=\mathrm{Im}d^{n-1}$.

If $\mathcal{C}$ is an abelian category and $\mathbf{f}:R\text{-}\mathrm{Mod}\rightarrow \mathcal{C}$ is an additive covariant functor, $\mathbf{f}(\mathbf{C})$ will be the complex
$$\cdots \longrightarrow \mathbf{f}(C^{n-1})\stackrel{\mathbf{f}(d^{n-1})}\longrightarrow
\mathbf{f}(C^{n})\stackrel{\mathbf{f}(d^{n})}\longrightarrow\mathbf{f}(C^{n+1})\stackrel{}\longrightarrow\cdots.$$ We say that $\mathbf{C}$ is $\mathbf{f}$-exact if $\mathbf{f}(\mathbf{C})$ is exact.

Let $\mathcal{C}$ be a Grothendieck category. $\mathcal{C}$ is said to be Gorenstein if it satisfies:

(1) The classes of all objects with finite projective dimension and with finite injective dimension coincide.

(2) The finitistic projective and injective dimensions of $\mathcal{C}$ are finite.

(3) $\mathcal{C}$ has a generator with finite projective dimension.

\begin{definition}([5, Definition 2.1]) A ring $R$ is said to be left Gorenstein regular if the category $R$-Mod is Gorenstein.
\end{definition}
Each Gorenstein ring (that is, a two sided noetherian ring with finite left and right self-injective dimensions) is left and right Gorenstein regular (see [7, Theorem 9.1.11]), and the converse is true precisely when the ring is two sided noetherian. A equivalent formulation for Gorenstein regular rings is given in [5, Proposition 2.2], which is more convenient to use. That is, a ring $R$ is left (resp. right) Gorenstein regular if and only if each projective left (resp. right) $R$-module has finite injective dimension and each injective left (resp. right) $R$-module has finite projective dimension.

Let $A$ and $B$ be rings and $U$ a $(B,A)$-bimodule. We denote by $T$ the triangular matrix ring $\left(\begin{smallmatrix}  A & 0 \\  U & B \\\end{smallmatrix}\right)$. According to [13, Theorem 1.5] that $T$-Mod is equivalent to the category whose objects are triples $M=\binom{M_{1}}{M_{2}}_{\varphi^{M}}$, where $M_{1}\in A$-$\mathrm{Mod}$, $M_{2}\in B$-$\mathrm{Mod}$ and $\varphi^{M}:U\otimes M_{1}\rightarrow M_{2}$ is a $B$-homomorphism, and whose morphisms between two objects $\binom{M_{1}}{M_{2}}_{\varphi^{M}}$ and $N=\binom{N_{1}}{N_{2}}_{\varphi^{N}}$~are pairs $\binom{f_{1}}{f_{2}}$ such that $f_{1}\in \mathrm{Hom}_{A}(M_{1},N_{1})$, $f_{2}\in \mathrm{Hom}_{B}(M_{2},N_{2})$, satisfying that the diagram
$$\xymatrix{
  U\otimes M_{1} \ar[d]_{\varphi^{M}} \ar[r]^{1_{U}\otimes f_{1}} & U\otimes N_{1} \ar[d]^{ \varphi^{N}} \\
  M_{2} \ar[r]^{f_{2}} & N_{2}   }
$$
is commutative. In the rest of the paper we identify $T$-Mod with this category and, whenever there is no possible confusion, we omit the homomorphism $\varphi^{M}$. Consequently, through the paper, a left $T$-module is a pair $\binom{M_{1}}{M_{2}}$. Given such a module $M$, we denote by $\widetilde{\varphi^{M}}$ the morphism from $M_{1}$ to $\mathrm{Hom}_{B}(U,M_{2})$ given by $\widetilde{\varphi^{M}}(m)(u)=\varphi^{M}( u\otimes m)$ for each $m\in M_{1},~u\in U$.

There are some functors between the category $T$-module and the product $A$-$\mathrm{Mod}\times B$-$\mathrm{Mod}$:

$\bullet$ $\mathbf{p}:A$-$\mathrm{Mod}\times B$-$\mathrm{Mod}\rightarrow T$-$\mathrm{Mod}$ is defined as follows: for each object
$(M_{1},M_{2})$ of $A$-$\mathrm{Mod}\times B$-$\mathrm{Mod}$, let $\mathbf{p}(M_{1},M_{2})=~\binom{M_{1}}{(U\otimes M_{1})\oplus M_{2}}$, with the obvious map. And, for any morphism $(f_{1},f_{2})$ in $A$-$\mathrm{Mod}\times B$-$\mathrm{Mod}$, let $\mathbf{p}(f_{1},f_{2})=~\binom{f_{1}}{(U\otimes f_{1})\oplus f_{2}}$.

$\bullet$ $\mathbf{h}:A$-$\mathrm{Mod}\times B$-$\mathrm{Mod}\rightarrow T$-$\mathrm{Mod}$
is defined as follows: for each object $(M_{1},M_{2})$ of $A$-$\mathrm{Mod}\times B$-$\mathrm{Mod}$, let $\mathbf{h}(M_{1},M_{2})=~\binom{M_{1}\oplus\mathrm{Hom}_{B}(U,M_{2})}{ M_{2}}$ with the obvious map.  And, for any morphism $(f_{1},f_{2})$ in $A$-$\mathrm{Mod}\times B$-$\mathrm{Mod}$, let $\mathbf{h}(f_{1},f_{2})=~\binom{f_{1}\oplus\mathrm{Hom}_{B}(U,f_{2})}{ f_{2}}$.

$\bullet$ $\mathbf{q}:T$-$\mathrm{Mod} \rightarrow A$-$\mathrm{Mod}\times B$-$\mathrm{Mod}$ is defined, for each left $T$-module $\binom{M_{1}}{M_{2}}$, as $\mathbf{q}\binom{M_{1}}{M_{2}}=(M_{1},M_{2})$, and for any morphism $\binom{f_{1}}{f_{2}}$ in $T$-Mod as $\mathbf{q}\binom{f_{1}}{f_{2}}=(f_{1},f_{2})$.

It is easy to see that $\mathbf{p}$ is a left adjoint of $\mathbf{q}$, $\mathbf{h}$ is a right adjoint of $\mathbf{q}$, and that $\mathbf{q}$ is exact. In particular, $\mathbf{p}$ preserves projective
 objects and $\mathbf{h}$ preserves injective objects.

Note that a sequence of $T$-modules
$$0\rightarrow\binom{M_{1}'}{M_{2}'}\rightarrow\binom{M_{1}}{M_{2}}\rightarrow\binom{M_{1}''}{M_{2}''}\rightarrow0$$
is exact if and only if both sequence $0\rightarrow M_{1}'\rightarrow M_{1}\rightarrow M_{1}''\rightarrow0$ of $A$-modules
and $0\rightarrow M_{2}'\rightarrow M_{2}\rightarrow M_{2}''\rightarrow0$ of $B$-modules are exact.
\section{Gorenstein projective (resp. injective) dimensions }
~~In this section, we describe explicitly the structure of left
$T$-modules of finite Gorenstein projective (resp.
injective) dimensions. We start with the following lemma, which is useful in the following arguments.
\begin{lemma}([5, Theorem 3.5]) Let $U_{A}$ have finite flat dimension, $_{B}U$ have finite projective dimension and $B$ be left Gorenstein regular. Let $M=\binom{M_{1}}{M_{2}}_{\varphi^{M}}$  be a left $T$-module. Then the following are equivalent:

$(1)~M$ is Gorenstein projective;

$(2)~M_{1}$ and $\frac{M_{2}}{\mathrm{Im}\varphi^{M}}$
are Gorenstein projective and the $\varphi^{M}$ is injective.
\end{lemma}
\begin{proposition} Let $B$ be a left Gorenstein regular ring, $U$ a $(B,A)$-bimodule, $M_{1}\in\mathcal{GP}(A)$. If $U$ is projective as left $B$-module and has finite flat dimension as right $A$-module, then $U\otimes_{A} M_{1}\in\mathcal{GP}(B).$
\end{proposition}
\noindent{\bf{proof.}} Assume that $U$ is projective as left $B$-module and has finite flat dimension as right $A$-module, $M_{1}\in\mathcal{GP}(A)$, let
$$\mathbf{P}:~~~~\cdots\longrightarrow P^{-1}\longrightarrow P^{0}\stackrel{\partial^{0}}\longrightarrow P^{1}\longrightarrow\cdots ~~~$$
be an exact sequence consisting of projective left $A$-modules which is $\mathrm{Hom}_{B}(-,Q)$-exact for each projective left $A$-module and such that $M_{1}\cong\mathrm{Ker}\partial^{0}_{\mathbf{P}}$. Using the hypothesis and [5, Lemma 2.3], we get that $U\otimes_{A} \mathbf{P}$ is an exact sequence of projective $B$-module such that $U\otimes M_{1}\cong\mathrm{Ker}(1_{U}\otimes\partial^{0}_{\mathbf{P}})$. Since $B$ is left Gorenstein regular, we know that $U\otimes_{A} M_{1}\in\mathcal{GP}(B)$ by [5, Proposition 2.6].\hfill$\Box$\\

In the following we consider that if left $T$-module $M=\binom{M_{1}}{M_{2}}_{\varphi^{M}}$ is Gorenstein projective, whether $M_{2}$ is a Gorenstein projective $B$-module.
\begin{proposition} Let $T$ be left Gorenstein regular, $U$ be projective as left $B$-module and have finite flat dimension as right $A$-module. If $G=\binom{G_{1}}{G_{2}}_{\varphi^{G}}\in\mathcal{GP}(T)$, then $G_{1}\in\mathcal{GP}(A)$, $G_{2}\in\mathcal{GP}(B)$ and the morphism $\varphi^{G}$ is injective.
\end{proposition}
\noindent{\bf{PROOF.}} Suppose $G=\binom{G_{1}}{G_{2}}_{\varphi^{G}}\in\mathcal{GP}(T)$. By Lemma 2.1, there exists a short exact sequence
$$0 \longrightarrow
U\otimes G_{1}\stackrel{\varphi^{G}}\longrightarrow G_{2}\longrightarrow \frac{G_{2}}{\mathrm{Im}\varphi^{G}}\stackrel{}\longrightarrow
0\\,$$
where $G_{1}$ and $\frac{G_{2}}{\mathrm{Im}\varphi^{G}}$ are Gorenstein projective. Since $U$ is projective as left $B$-module and has finite flat dimension as right $A$-module, it follows from Proposition 2.2 that $U\otimes_{A} G_{1}\in\mathcal{GP}(B)$. Thus $G_{2}\in\mathcal{GP}(B)$ by [17, Theorem 2.5].\hfill$\Box$\\

Asadollahi and Salarian [1]
establish a relationship between the projective (resp. injective) dimension
of modules over $T$ and over $A$ and $B$. Let $n$ be a non-negative
integer. They proved that $\mathrm{pd}\binom{M_{1}}{M_{2}}_{\varphi^{M}}\leq n$ if and
only if $\mathrm{pd}(M_{1})\leq n$,
$\mathrm{pd}(\frac{M_{2}}{\mathrm{Im}\varphi^{M}})\leq n$ and the map related
to the $n$-th syzygy of $\binom{M_{1}}{M_{2}}_{\varphi^{M}}$ is
injective. We have similar arguments about Gorenstein projective
dimension and Gorenstein injective dimension.\\

The following lemma is quoted from [1].
\begin{lemma}Let $J=Te_{1}$, where $e_{1}=\left( \begin{smallmatrix}  1 & 0 \\ 0 & 0 \\  \end{smallmatrix}\right)$. Then for each left
 $T$-module $X=\binom{X_{1}}{X_{2}}_{\varphi^{X}}$, there is a isomorphism $\frac{T}{J}\otimes X\cong \binom{0}{\frac{X_{2}}{\mathrm{Im}\varphi^{X}}}$.
\end{lemma}
\begin{theorem} Let $n$ be a non-negative integer, $B$ a left Gorenstein regular ring, $U_{A}$ have finite flat dimension and $_{B}U$ have
finite projective dimension. Let
$M=\binom{M_{1}}{M_{2}}_{\varphi^{M}}$  be a left $T$-module. Then
$\mathrm{Gpd}(M)\leq n$ if and only if $\mathrm{Gpd}(M_{1})\leq n$,
$\mathrm{Gpd}(\frac{M_{2}}{\mathrm{Im}\varphi^{M}})\leq n$ and if
$K=\binom{K_{1}}{K_{2}}_{\varphi^{K}}$ is a n-th syzygy of $M$, then
$\varphi^{K}$ is injective.
\end{theorem}
\noindent{\bf{PROOF.}} $(\Rightarrow)$ Let $\mathrm{Gpd}(M)\leq n$. Then there exists an exact sequence of
$T$-modules
$$\mathbf{P}:~0 \longrightarrow
\binom{G^{n}_{1}}{G^{n}_{2}}\stackrel{\binom{\partial^{n}_{1}}{\partial^{n}_{2}}}\longrightarrow \cdots\longrightarrow\binom{G^{1}_{1}}{G^{1}_{2}}\stackrel{\binom{\partial^{1}_{1}}{\partial^{1}_{2}}}\longrightarrow \binom{G^{0}_{1}}{G^{0}_{2}}\stackrel{\binom{\partial^{0}_{1}}{\partial^{0}_{2}}}\longrightarrow
\binom{M_{1}}{M_{2}}\longrightarrow 0,~~~~~~~~~~~~~~~~~~~\\$$
where $\binom{G^{i}_{1}}{G^{i}_{2}}\in\mathcal{GP}(T)$ for $i=0,\cdots,n$. Thus the sequence
$$\mathbf{P}_{1}: 0 \longrightarrow
G^{n}_{1}\longrightarrow \cdots\longrightarrow
G^{1}_{1}\stackrel{}\longrightarrow M_{1}\longrightarrow
0~~~~~~~~~~~~~~~~~~~~~~~~~~~~~~~~~~~~~~~~~~~~~~~~\\$$
is exact. Since
$\binom{G^{i}_{1}}{G^{i}_{2}}\in\mathcal{GP}(T)$, then
$G^{i}_{1}\in\mathcal{GP}(A)$ by Lemma 2.1, and so $\mathrm{Gpd}(M_{1})\leq n$. By Lemma 2.4, because $\frac{T}{J}$ is projective in $T$-Mod, we
 can apply  the functor
$\frac{T}{J}\otimes-$ on $\mathbf{P}$ to get the following exact
sequence
$$\overline{\mathbf{P}_{2}}:~~~~ 0 \longrightarrow
\frac{G^{n}_{2}}{\mathrm{Im}\varphi^{n}}\longrightarrow \cdots\longrightarrow \frac{G^{0}_{2}}{\mathrm{Im}\varphi^{0}}\stackrel{}\longrightarrow \frac{M_{2}}{\mathrm{Im}\varphi^{M}}\longrightarrow
0,~~~~~~~~~~~~~~~~~~~~~~~~~~~~~~\\$$
where $\frac{G^{i}_{2}}{\mathrm{Im}\varphi^{i}}\in\mathcal{GP}(B)$ for $i=0,\cdots,n$. Thus $\mathrm{Gpd}(\frac{M_{2}}{\mathrm{Im}\varphi^{M}})\leq n$.

If $K=\binom{K_{1}}{K_{2}}_{\varphi^{K}}$ is a $n$-th syzygy of $M$, since $\mathrm{Gpd}(M)\leq n$, $\binom{K_{1}}{K_{2}}$  has to be Gorenstein projective by [17, Theorem 2.20], and the
 morphism $\varphi^{K}$ is injective by Lemma 2.1.

$(\Leftarrow)$ Let $\mathrm{Gpd}(M_{1})\leq n$, $\mathrm{Gpd}(\frac{M_{2}}{\mathrm{Im}\varphi^{M}})\leq n$, and $\varphi^{K}$ be injective for any $n$-th syzygy $\binom{K_{1}}{K_{2}}_{\varphi^{K}}$ of $M$.
 Then there exists an exact sequence
$$\mathbf{P}':~~~~~0 \longrightarrow
\binom{K^{n}_{1}}{K^{n}_{2}}\stackrel{}\longrightarrow \binom{P^{n-1}_{1}}{P^{n-1}_{2}}\longrightarrow\cdots\longrightarrow\binom{P^{1}_{1}}{P^{1}_{2}}\stackrel{}\longrightarrow \binom{P^{0}_{1}}{P^{0}_{2}}\stackrel{}\longrightarrow
\binom{M_{1}}{M_{2}}\longrightarrow 0\\.$$
Since every projective $T$-module is Gorenstein projective, it suffices to verify  $\binom{K^{n}_{1}}{K^{n}_{2}}\in\mathcal{GP}(T)$.

The sequence above induces the following exact sequence
$$ 0 \longrightarrow
K^{n}_{1}\longrightarrow P^{n-1}_{1}\longrightarrow\cdots\longrightarrow P^{0}_{1}\stackrel{}\longrightarrow M_{1}\longrightarrow
0.\\$$
As $\mathrm{Gpd}(M_{1})\leq n$, we know that $K^{n}_{1}\in\mathcal{GP}(A)$.
Apply the functor $\frac{T}{J}\otimes-$ on $\mathbf{P}'$, we get the following exact sequence
$$~~~~ 0 \longrightarrow
\frac{K^{n}_{2}}{\mathrm{Im}\varphi^{K}}\longrightarrow \cdots\longrightarrow \frac{P^{0}_{2}}{\mathrm{Im}\varphi^{0}}\stackrel{}\longrightarrow \frac{M_{2}}{\mathrm{Im}\varphi^{M}}\longrightarrow
0.~~~~~~~~~~\\$$
We have $\frac{K^{n}_{2}}{\mathrm{Im}\varphi^{K}}\in\mathcal{GP}(B)$ since  $\mathrm{Gpd}(\frac{M_{2}}{\mathrm{Im}\varphi^{M}})\leq n$. We get
 that $\binom{K^{n}_{1}}{K^{n}_{2}}\in\mathcal{GP}(T)$ by Lemma 2.1. Thus $\mathrm{Gpd}(M)\leq n$.\hfill$\Box$\\

Using the dual argument of Lemma 2.1 and Theorem 2.5,
one can prove the following theorem.
\begin{theorem} Let $n$ be a non-negative integer, $A$ a left Gorenstein regular ring, $U_{A}$ have finite flat dimension and $_{B}U$ have finite projective dimension. Let $M=\binom{M_{1}}{M_{2}}_{\varphi^{M}}$  be a left $T$-module. Then $\mathrm{Gid}(M)\leq n$ if and only if  $\mathrm{Gid}(M_{2})\leq n$,~$\mathrm{Gid}(\mathrm{Ker}\widetilde{\varphi^{M}})\leq n$ and if $\binom{L_{1}}{L_{2}}_{\varphi^{L}}$ is a n-th cosyzygy of $M$, then $\widetilde{\varphi^{L}}$ is surjective.
\end{theorem}

For finiteness of Gorenstein projective (resp. injective) dimensions, we have the following result.

\begin{theorem} Let $U_{A}$ have finite flat dimension and $_{B}U$ have finite projective dimension. Let $M=\binom{M_{1}}{M_{2}}_{\varphi^{M}}$  be a left $T$-module.

(1) Suppose that $B$ is left Gorenstein regular, then $\mathrm{Gpd}(M)< \infty$ if and only if $\mathrm{Gpd}(M_{1})< \infty$,~$\mathrm{Gpd}(M_{2})< \infty.$

(2) Suppose that $A$ is left Gorenstein regular, then
$\mathrm{Gid}(M)< \infty$ if and only if $\mathrm{Gid}(M_{1})<
\infty$,~$\mathrm{Gid}(M_{2})< \infty.$
\end{theorem}
\noindent{\bf{PROOF.}} (1)$(\Rightarrow)$ Suppose $\mathrm{Gpd}(M)< \infty$. By Theorem 2.5, we get $\mathrm{Gpd}(M_{1})\leq n$. As a consequence of [9, Theorem 2.27], $B$ is left Gorenstein regular if and only if $\mathrm{glGpd}(B)< \infty$ and $\mathrm{glGid}(B)< \infty$, Thus $\mathrm{Gpd}(M_{2})< \infty$
by the hypothesis.

$(\Leftarrow)$ Suppose that $\mathrm{Gpd}(M_{1})< \infty$
and $\mathrm{Gpd}(M_{2})< \infty$. Fix a Gorenstein
projective resolution of $M_{1}$,
$$~~~~ 0 \longrightarrow
G^{n}_{1}\stackrel{\partial^{n}_{1}}\longrightarrow \cdots\longrightarrow G^{0}_{1}\stackrel{\partial^{0}_{1}}\longrightarrow M_{1}\longrightarrow0~~~~~\\$$
and a Gorenstein projective presentation $G^{0}_{2}\stackrel{\partial^{0}_{2}}\rightarrow M_{2}$ of $M_{2}$. Then $\mathbf{p}(\partial^{0}_{1},\partial^{0}_{2})$ is a Gorenstein projective presentation of $M$ in $T$-Mod and, if  $K^{0}=\binom{K^{0}_{1}}{K^{0}_{2}}$ is its kernel, there exists short exact sequences
$$~~~~ 0 \longrightarrow
K^{0}_{1}\longrightarrow  G_{1}\longrightarrow M_{1}\longrightarrow 0~~~~~\\$$
and
$$~~~~ 0 \longrightarrow
K^{0}_{2}\longrightarrow  G_{2}\oplus(U\otimes_{A}G_{1})\longrightarrow M_{2}\longrightarrow 0.~~~~~\\$$
Then $K^{0}_{1}$ and $K^{0}_{2}$ have finite Gorenstein projective dimension by [17, Theorem 2.24]. Now take $\partial^{1}_{2}:G^{1}_{2}\rightarrow K^{0}_{2}$ a projective presentation. Since $K^{0}_{1}$ is the kernel of $\partial^{0}_{1}$, $\mathbf{p}(\partial^{1}_{1},\partial^{1}_{2})$ is a projective presentation of $K^{0}$ and, reasoning as before, its kernel, say $K^{1}=\binom{K^{1}_{1}}{K^{1}_{2}}$, satisfy that $K^{1}_{1}$ and $K^{1}_{2}$ have finite Gorenstein projective dimension. Repeating this procedure, we construct a Gorenstein projective resolution of $M$
$$~~~~\cdots \longrightarrow
\mathbf{p}(G^{1}_{1},G^{1}_{2})\longrightarrow  \mathbf{p}(G^{0}_{1},G^{0}_{2})\longrightarrow M\longrightarrow 0,~~~~~\\$$
such that $K^{m}=\binom{K^{m}_{1}}{K^{m}_{2}}$ is the kernel of $\mathbf{p}(\partial^{m}_{1},\partial^{m}_{2})$ for each $m\in\mathbb{N}$, both $K^{m}_{1}$ and $K^{m}_{2}$ have finite Gorenstein dimensions. Since $G^{n+1}_{1}=0$, $\mathrm{Ker}\mathbf{p}(\partial^{n+1}_{1},\partial^{n+1}_{2})=\binom{0}{K^{n+1}_{2}}$. As $K^{n+1}_{2}$ has finite Gorenstein projective dimension in $B$-Mod, we have the following exact sequence
$$~~~~ 0 \longrightarrow
Q^{n+m}_{2}\longrightarrow \cdots\longrightarrow Q^{n+2}_{2}\longrightarrow K^{n+1}_{2}\longrightarrow0,~~~~~\\$$
which induces the finite resolution in $T$-Mod
$$~~~~0 \longrightarrow
\mathbf{p}(0,Q^{n+m}_{2})\longrightarrow\cdots\longrightarrow\mathbf{p}(0,Q^{n+2}_{2})\longrightarrow \binom{0}{K^{n+1}_{2}}\longrightarrow 0.~~~~~\\$$
This means that $\binom{0}{K^{n+1}_{2}}$ has finite Gorenstein projective dimension, which implies that $\mathrm{Gpd}(M)< \infty$.

The proof of (2) is similar to the way we do in (1), we
can prove $\mathrm{Gid}(M)< \infty$ if and only if
$\mathrm{Gid}(M_{1})< \infty$,~$\mathrm{Gid}(M_{2})<
\infty$.\hfill$\Box$\\

It follows from [2, Lemma 3.4] that if $_{B}U$ is projective, $\mathrm{pd}(M_{1})\leq n$ and $\mathrm{pd}(M_{2})\leq n$, then $\mathrm{pd}\binom{M_{1}}{M_{2}}\leq n+1$. Moreover, we have the following result, which is a more explicit version of Theorem 2.7 under some additional conditions.
\begin{lemma}Let $_{B}U$ have finite projective dimension, $U_{A}$ have finite flat dimension, $B$ be a left Gorenstein
regular ring. Let $M=\binom{M_{1}}{M_{2}}_{\varphi^{M}}\in\mathcal{GP}(T)$, $U\otimes _{A}M_{1}$ have finite
projective dimension. Then there exists $P\in\mathcal{GP}(A)$ and $Q\in\mathcal{GP}(B)$
such that $M=\mathbf{p}(P,Q)$.
\end{lemma}
\noindent{\bf{PROOF.}} There exists an exact sequence
$$0 \longrightarrow
U\otimes M_{1}\stackrel{\varphi^{M}}\longrightarrow
M_{2}\longrightarrow
\frac{M_{2}}{\mathrm{Im}\varphi^{M}}\stackrel{}\longrightarrow 0\\$$ in which
$M_{1}$ and $\frac{M_{2}}{\mathrm{Im}\varphi^{M}}$ are Gorenstein projective.
Since  $U\otimes _{A}M_{1}$ has finite
projective dimension, we know by [17, Theorem 2.20] that
$\mathrm{Ext}^{1}_{B}(\frac{M_{2}}{\mathrm{Im}\varphi^{M}},U\otimes
M_{1})=0$, which means the above sequence splits, and
$M_{2}=(U\otimes M_{1})\oplus \frac{M_{2}}{\mathrm{Im}\varphi^{M}}$. Let
$P=M_{1},~Q=\frac{M_{2}}{\mathrm{Im}\varphi^{M}}$, then
$$M=\binom{M_{1}}{M_{2}}=\binom{M_{1}}{U\otimes M_{1}}\oplus\binom{0}{\frac{M_{2}}{\mathrm{Im}\varphi^{M}}}=\binom{P}{U\otimes P}\oplus\binom{0}{Q}=\mathbf{p}(P,Q).$$~\hfill$\Box$\\

Dually, we have the following result for
Gorenstein injective $T$-modules.
\begin{lemma}Let $_{B}U$ have finite projective dimension, $U_{A}$ have finite flat dimension, $A$ be a left Gorenstein
regular ring. Let $M=\binom{M_{1}}{M_{2}}_{\varphi^{M}}\in\mathcal{GI}(T)$, $\mathrm{Hom}_{B}(U,M_{2})$ have finite injective dimension.
Then there exists $E\in\mathcal{GI}(A)$ and $I\in\mathcal{GP}(B)$
such that $M=\mathbf{h}(E,I)$.
\end{lemma}

\begin{proposition} Let $U_{A}$ has finite flat dimension, $B$ be a left Gorenstein regular, $M=\binom{M_{1}}{M_{2}}$ a $T$-module, $Gpd(M_{1})\leq n$, $Gpd(M_{2})\leq n$, $U\otimes _{A}M_{1}$ has finite
projective dimension. If
$_{B}U$ is projective, then $\mathrm{Gpd}(M)\leq n+1$.
\end{proposition}
\noindent{\bf{PROOF.}} By Lemma 2.8, a Gorenstein projective
resolution of the $T$-module $M$ can be written in the following
form
$$\cdots \rightarrow \mathbf{p}(P_{n},Q_{n}) \rightarrow\cdots\rightarrow\mathbf{p}(P_{1},Q_{1}) \rightarrow\mathbf{p}(P_{0},Q_{0}) \rightarrow M\rightarrow 0,$$
where $P_{i}\in\mathcal{GP}(A)$ and $Q_{i}\in\mathcal{GP}(B)$ for
$i>0$. Let $\binom{K_{1}}{K_{2}}$ be the n-th syzygy of $M$. We show
that $\binom{K_{1}}{K_{2}}$ has Gorenstein projective dimension at most one.
The above resolution induces the following exact sequence
$$0\rightarrow K_{1}\rightarrow P_{n-1}\rightarrow \cdots\rightarrow P_{1}\rightarrow P_{0}\rightarrow M_{1}\rightarrow0 $$
of $A$-modules and an exact sequence
$$0\rightarrow K_{2}\rightarrow (U\otimes P_{n-1})\oplus Q_{n-1}\rightarrow\cdots\rightarrow(U\otimes P_{1})\oplus Q_{1}\rightarrow(U\otimes P_{0})\oplus Q_{0}\rightarrow M_{2}\rightarrow0$$
of $B$-module. We know that $U\otimes P_{n-1}\in\mathcal{GP}(B)$ by
Proposition 2.2 for $i=0,\cdots,n-1$. The projectivity of $_{B}U$ in
conjunction with our assumption on the Gorenstein projective
dimensions of $M_{1}$ and $M_{2}$ imply that $K_{1}$and $K_{2}$ are
Gorenstein projective. Now consider the following short exact
sequence of $T$-modules
$$0\longrightarrow\binom{0}{U\otimes K_{1}}\longrightarrow\binom{K_{1}}{U\otimes K_{1}}\oplus\binom{0}{K_{2}}\longrightarrow\binom{K_{1}}{K_{2}}\longrightarrow0$$
Since $_{B}U$ is projective, we deduce that the first two left terms
are Gorenstein projective, and so $\mathrm{Gpd}\binom{K_{1}}{K_{2}}\leq 1$. Therefore $\mathrm{Gpd}(M)\leq n+1$.\hfill$\Box$\\
\section{Gorenstein flat dimensions }
In the following, we describe Gorenstein flat $T$-module over triangular
matrix ring. It has proved in [5, Theorem 3.1] that the triangular
matrix ring $T$ is left Gorenstein regular if and only if $A$ and $B$
are Gorenstein regular when $_{B}U$ has finite projective dimension
and $U_{A}$ has finite flat dimension. It have analogous versions for right Gorenstein
regular rings. That is, if $U_{A}$ has finite projective dimension and $_{B}U$ has finite flat dimension, then the triangular matrix ring $T$ is right Gorenstein regular if and only if $A$ and $B$ are right Gorenstein regular.
\begin{lemma}([12, Proposition 1.8]) Let $_{B}U$ ($U_{A}$) be finitely generated. Then
$T$ is left (right) noetherian if and only if $A$ and $B$ are left
(right) noetherian.
\end{lemma}
\begin{proposition} Let $U_{A}$ and $_{B}U$ be finitely generated and have finite projective dimensions. Then $T$ is Gorenstein if and only if $A$ and $B$ are Gorenstein.
\end{proposition}
\noindent{\bf{PROOF.}} We know that, by [5], $T$ is Gorenstein if and only if $T$ is left and right Gorenstein regular when $T$ is a two sided noetherian ring. Then we get that $T$ is Gorenstein if and only if $A$ and $B$ are Gorenstein by Lemma 3.1 and [5, Theorem 3.1].\hfill$\Box$
\begin{lemma} ([10, Proposition 1.14]) Let $F=\binom{F_{1}}{F_{2}}_{\varphi^{F}}$ be a left $T$-module. Then $F$ is flat if and only if $F_{1}$ and $\frac{F_{2}}{Im\varphi^{F}}$ are flat and morphism $\varphi^{F}$ is injective.
\end{lemma}
\begin{proposition}  Let $n$ be a non-negtive integer, $M=\binom{M_{1}}{M_{2}}_{\varphi^{M}}$  a left $T$-module. Then $\mathrm{fd}(M)\leq n$ if and only if $\mathrm{fd}(M_{1})\leq n$, $\mathrm{fd}(\frac{M_{2}}{\mathrm{Im}\varphi^{M}})\leq n$, and if $\binom{K_{1}}{K_{2}}_{\varphi^{K}}$ is a $n$-th syzygy of $M$, then $\varphi^{K}$ is injective.
\end{proposition}
\noindent{\bf{PROOF.}} $(\Rightarrow)$ Let $\mathrm{fd}(M)\leq n$.
There exists an exact sequence of $T$-modules
$$\mathbf{F}:~0 \longrightarrow
\binom{F^{n}_{1}}{F^{n}_{2}}\stackrel{\binom{\partial^{n}_{1}}{\partial^{n}_{2}}}\longrightarrow
\cdots\longrightarrow\binom{F^{1}_{1}}{F^{1}_{2}}\stackrel{\binom{\partial^{1}_{1}}{\partial^{1}_{2}}}\longrightarrow
\binom{F^{0}_{1}}{F^{0}_{2}}\stackrel{\binom{\partial^{0}_{1}}{\partial^{0}_{2}}}\longrightarrow
\binom{M_{1}}{M_{2}}\longrightarrow 0,~~~~~~~~~~~~~~~~~~~\\$$ where
$\binom{F^{i}_{1}}{F^{i}_{2}}$ is flat for $i=0,\cdots,n$, which induces the exact sequence
$$\mathbf{F}_{1}: 0 \longrightarrow
F^{n}_{1}\longrightarrow \cdots\longrightarrow
F^{1}_{1}\stackrel{}\longrightarrow M_{1}\longrightarrow
0.~~~~~~~~~~~~~~~~~~~~~~~~~~~~~~~~~~~~~~~~~~~~~~~~\\$$ Since
$\binom{F^{i}_{1}}{F^{i}_{2}}$ is flat, $F^{i}_{1}$ is flat by
Lemma 3.3, and so $\mathrm{fd}(M_{1})\leq n$. By Lemma 2.4, we apply the functor
$\frac{T}{J}\otimes-$ on $\mathbf{F}$ to get the following exact
sequence
$$\overline{\mathbf{F}_{2}}:~~~~ 0 \longrightarrow
\frac{F^{n}_{2}}{\mathrm{Im}\varphi^{n}}\longrightarrow \cdots\longrightarrow
\frac{F^{0}_{2}}{\mathrm{Im}\varphi^{0}}\stackrel{}\longrightarrow
\frac{M_{2}}{\mathrm{Im}\varphi^{M}}\longrightarrow
0~~~~~~~~~~~~~~~~~~~~~~~~~~~~~~\\$$
where $\frac{F^{i}_{2}}{\mathrm{Im}\varphi^{i}}$ is flat for $i=0,\cdots,n$. Thus $\mathrm{fd}(\frac{M_{2}}{\mathrm{Im}\varphi^{M}})\leq n$. If $\binom{K_{1}}{K_{2}}_{\varphi^{K}}$ is n-th syzygy of $M$, $\binom{K_{1}}{K_{2}}_{\varphi^{K}}$ is flat since $\mathrm{fd}(M)\leq n$, and morphism $\varphi^{K}$ is injective by Lemma 3.3.

$(\Leftarrow)$ Suppose $\mathrm{fd}(M_{1})\leq n$,
$\mathrm{fd}(\frac{M_{2}}{\mathrm{Im}\varphi^{M}})\leq n$, and
$\binom{K_{1}}{K_{2}}_{\varphi^{K}}$ is n-th syzygy
of $M$, $\varphi^{K}$ is injective. Then there exists an exact
sequence
$$\mathbf{F}':~~~~~0 \longrightarrow
\binom{K^{n}_{1}}{K^{n}_{2}}\stackrel{}\longrightarrow
\binom{F^{n-1}_{1}}{F^{n-1}_{2}}\longrightarrow\cdots\longrightarrow\binom{F^{1}_{1}}{F^{1}_{2}}\stackrel{}\longrightarrow
\binom{F^{0}_{1}}{F^{0}_{2}}\stackrel{}\longrightarrow
\binom{M_{1}}{M_{2}}\longrightarrow 0\\.$$ Then it only need to
verify that $\binom{K^{n}_{1}}{K^{n}_{2}}$ is flat. We have an exact
sequence
$$ 0 \longrightarrow
K^{n}_{1}\longrightarrow
F^{n-1}_{1}\longrightarrow\cdots\longrightarrow
F^{0}_{1}\stackrel{}\longrightarrow M_{1}\longrightarrow 0.\\$$ As
$\mathrm{fd}(M_{1})\leq n$, we know that $K^{n}_{1}$ is flat. Apply
the functor $\frac{T}{J}\otimes-$ on $\mathbf{F}'$, by Lemma 2.4 we get the following
exact sequence
$$~~~~ 0 \longrightarrow
\frac{K^{n}_{2}}{\mathrm{Im}\varphi^{K}}\longrightarrow \cdots\longrightarrow
\frac{P^{0}_{2}}{\mathrm{Im}\varphi^{0}}\stackrel{}\longrightarrow
\frac{M_{2}}{\mathrm{Im}\varphi^{M}}\longrightarrow 0,~~~~~~~~~~\\$$
we obtain that $\frac{K^{n}_{2}}{\mathrm{Im}\varphi^{K}}$ is flat since  $\mathrm{fd}(\frac{M_{2}}{\mathrm{Im}\varphi^{M}})\leq n$. The previous Lemma implies that $\binom{K^{n}_{1}}{K^{n}_{2}}$ is flat. Thus $\mathrm{fd}(M)\leq n$.\hfill$\Box$
\begin{proposition}Let $M=\binom{M_{1}}{M_{2}}_{\varphi^{M}}$ be a left $T$-module. Suppose that $_{B}U$ has finite flat dimension, Then $M$ has finite flat dimension if and only if $M_{1}$ and $M_{2}$ have finite flat dimension.
\end{proposition}
\noindent{\bf{PROOF.}} It is ture by combining Proposition 3.4 with [5, Proposition 2.8(1)].\hfill$\Box$

\begin{proposition} Let $U_{A}$ and $_{B}U$ have finite projective dimension, $T$ be a right Gorenstein regular ring, $M=\binom{M_{1}}{M_{2}}_{\varphi^{M}}$ a left $T$-module. If $M$ is $Gorenstein$ flat, then $M_{1}$ and $\frac{M_{2}}{\mathrm{Im}\varphi^{M}}$ are $Gorenstein$ flat and the morphism $\varphi^{M}$ is injective.
\end{proposition}
\noindent{\bf{PROOF.}} Suppose $M$ is Gorenstein flat in $T$-Mod
and let
$$~\mathbf{F}:~\cdots \longrightarrow
\binom{F^{n-1}_{1}}{F^{n-1}_{2}}\stackrel{\binom{\partial^{n-1}_{1}}{\partial^{n-1}_{2}}}
\longrightarrow\binom{F^{n}_{1}}{F^{n}_{2}}\stackrel{\binom{\partial^{n}_{1}}{\partial^{n}_{2}}}\longrightarrow
\binom{F^{n+1}_{1}}{F^{n+1}_{2}}\stackrel{}\longrightarrow
\cdots$$
be an exact sequence consisting of flat left $T$-modules
which is $E\otimes-$-exact for each injective right $T$-module $E$ and
such that $\mathrm{Ker}\binom{\partial^{0}_{1}}{\partial^{0}_{2}}\cong M$. Then
$$\mathbf{F}_{1}:~~~~ \cdots \longrightarrow
F^{n-1}_{1}\longrightarrow F^{n}_{1}\stackrel{}\longrightarrow
F_{n+1}\longrightarrow\cdots~~~~~~~~~~\\$$
is an exact sequence consisting of flat left $A$-modules with $\mathrm{Ker}\partial^{0}_{1}\cong M_{1}$. Moreover, $A$ is right Gorenstein regular, then for each injective right $A$-module $I$, $I$ has finite flat dimension, we get that $I\otimes\mathbf{F}_{1}$ is exact by [5, Lemma 2.3]. This means that $M_{1}$ is a Gorenstein flat left $A$-module.

Now note that for every morphism
$\binom{\partial^{n}_{1}}{\partial^{n}_{2}}:\binom{F^{n}_{1}}{F^{n}_{2}}\stackrel{}\longrightarrow
\binom{F^{n+1}_{1}}{F^{n+1}_{2}}$, we can construct the following
commutative diagram:
$$\CD
  0 @> >>U\otimes F^{n}_{1} @>\varphi^{n}>>F^{n}_{2}@>\pi^{n}>> \frac{F^{n}_{2}}{\mathrm{Im}\varphi^{n}}@> >> 0 \\
  @. @V 1_{U}\otimes\partial^{n}_{1} VV @V  \partial^{n}_{2} VV @V \overline{\partial^{n}_{2}} VV @.   \\
  0 @> >>U\otimes F^{n+1}_{1} @>\varphi^{n+1}>>F^{n+1}_{2}@>\pi^{n+1}>> \frac{F^{n+1}_{2}}{\mathrm{Im}\varphi^{n+1}}@> >> 0 ,
\endCD $$
where $\pi^{n}$ and $\pi^{n+1}$ are the canonical projections. Using
this fact, the complex $\mathbf{F}$ induce the complex
$$\overline{\mathbf{F}_{2}}:~~~~ \cdots \longrightarrow
\frac{F^{n-1}_{2}}{\mathrm{Im}\varphi^{n-1}}\stackrel{\overline{\partial^{n-1}_{2}}}\longrightarrow
\frac{F^{n}_{2}}{\mathrm{Im}\varphi^{n}}\stackrel{\overline{\partial^{n}_{2}}}\longrightarrow
\frac{F^{n+1}_{2}}{\mathrm{Im}\varphi^{n+1}}\longrightarrow\cdots,~~~~~~~~~~\\$$
where $\varphi^{i}$ is the structural map of the $T$-module
$\binom{F^{i}_{1}}{F^{i}_{2}}$ for each $i\in \mathbb{Z}$. We get that each $\frac{F^{n}_{2}}{\mathrm{Im}\varphi^{n}}$ is Gorenstein flat in
$B$-Mod, since $\binom{F^{i}_{1}}{F^{i}_{2}}$ is a Gorenstein flat
left $T$-module for each $i\in \mathbb{Z}$. Moreover, the complex
$\overline{\mathbf{F}_{2}}$ is exact, since there exists a short
exact sequence of complexes
$$~~~~0 \longrightarrow
U\otimes\mathbf{F}_{1}\longrightarrow
\mathbf{F}_{2}\stackrel{}\longrightarrow
\overline{\mathbf{F}_{2}}\longrightarrow0,~~~~~\\$$
in which $U\otimes\mathbf{F}_{1}$ is exact by [5, Lemma 2.3], and $\mathbf{F}_{2}$ is exact. It is easy to see that $\mathrm{Ker}\overline{\partial^{0}_{2}}\cong \frac{M_{2}}{Im\varphi^{M}}$. Since $B$ is right Gorenstein regular, for each injective right $B$-module $E$, $E$ has finite flat dimensions, then we get that $\overline{\mathbf{F}_{2}}$ is $E\otimes-$-exact by [5, Lemma 2.3]. Thus $\frac{M_{2}}{\mathrm{Im}\varphi^{M}}$ is Gorenstein flat.

Finally, we prove that morphism $\varphi^{M}$ is
injective. By [5, Lemma 2.3], $U\otimes\mathbf{F}_{1}$ is an exact
sequence. This means that if $\iota_{1}:M_{1}\rightarrow F^{0}_{1}$
is the inclusion, $1_{U}\otimes\iota_{1}$ is injective. But, since
$\binom{M_{1}}{M_{2}}_{\varphi^{M}}$ is a submodule of
$\binom{F^{0}_{1}}{F^{0}_{2}}_{\varphi^{F}}$, the following diagram
commutes:
$$\CD
 U\otimes M_{1} @>1_{U}\otimes\iota_{1}>>U\otimes F^{0}_{1}  \\
@V\varphi^{M} VV @V \varphi^{0} VV   \\
  M_{2} @>\iota_{2}>> F^{0}_{1} ,
\endCD $$
where $\iota_{2}$ is the inclusion. Since $\varphi^{0}$ is injective, and $\binom{F^{0}_{1}}{F^{0}_{2}}_{\varphi^{F}}$ is Gorenstein flat, we conclude that $\varphi^{M}$ is injective.\hfill$\Box$

In the following we will show that the converse of Proposition 3.6 is ture under some additional conditions.
\begin{lemma}([12]) Let $M=\binom{M_{1}}{M_{2}}_{\varphi^{M}}$ be a left $T$-module. Then $M$ is finitely generated if and only if $M_{1}$ and $\frac{M_{2}}{\mathrm{Im}\varphi^{M}}$ are finitely generated.
\end{lemma}
\begin{lemma} Let $U_{A}$ and $_{B}U$ be finitely generated and have finite projective dimension, $T$ be a Gorenstein ring, $F=\binom{F_{1}}{F_{2}}_{\varphi^{F}}$ a left $T$-module. If  $M_{1}$ and $M_{2}$ are Gorenstein flat left $A$-module and left $B$-module respectively, then $\mathbf{p}(M_{1},M_{2})$ is Gorenstein flat left $T$-module.
\end{lemma}
\noindent{\bf{PROOF.}} We first note that $A$ and $B$ are Gorenstein by Proposition 3.2. Then it follows from [6, Theorem 10.3.8] that both $M_{1}$ and $M_{2}$ are direct limit of finitely generated Gorenstein projective modules. But if $P_{1}$ is a finitely generated Gorenstein projective $A$-module, then $\binom{P_{1}}{U\otimes P_{1}}$ is a finitely generated Gorenstein projective $T$-modules. Suppose that $M_{1}=\underrightarrow{\mathrm{lim}}P^{i}_{1}$, where each module $P^{i}_{1}$, $i\in \mathbb{Z}$ is finitely generated Gorenstein projective, we get that $\binom{M_{1}}{U\otimes M_{1}}=\underrightarrow{\mathrm{lim}}\binom{P^{i}_{1}}{U\otimes P^{i}_{1}}$. Hence $\binom{M_{1}}{U\otimes M_{1}}\in\mathcal{GF}(T)$, reasoning as before, so is $\binom{0}{M_{2}}$. Then $\mathbf{p}(M_{1},M_{2})=\binom{M_{1}}{U\otimes M_{1}}\oplus\binom{0}{M_{2}}$ is Gorenstein flat left $T$-module.\hfill$\Box$
\begin{theorem} Let $U_{A}$ and $_{B}U$ be finitely generated and have finite projective dimension, $T$ be a Gorenstein ring, $F=\binom{F_{1}}{F_{2}}_{\varphi^{F}}$ a left $T$-module. Then $F\in\mathcal{GF}(T)$ if and only if $F_{1}\in\mathcal{GF}(A)$, $\frac{F_{2}}{\mathrm{Im}\varphi^{F}}\in\mathcal{GF}(B)$ and $\varphi^{F}$ is injective.
\end{theorem}
\noindent{\bf{PROOF.}} We see that if
$F=\binom{F_{1}}{F_{2}}$ is Gorenstein flat, then as a limit of
finitely generated Gorenstein projective T-modules
$P=\binom{P_{1}}{P_{2}}_{\varphi^{P}}$ (with $P_{1}$ necessarily a
finitely generated Gorenstein projective modules) we get that
$F_{1}$ is Gorenstein flat. We also note that each such $P$
satisfies that $\frac{P_{2}}{\mathrm{Im}\varphi^{P}}$ is finitely generated Gorenstein projective and $\varphi^{F}$ is injective. Since these
conditions are preserved by direct limits we get that $\frac{F_{2}}{\mathrm{Im}\varphi^{F}}\in\mathcal{GF}(B)$ and $\varphi^{F}$ is injective.

Conversely, we assume that
$F=\binom{F_{1}}{F_{2}}_{\varphi^{F}}$ satisfies that $F_{1}$ and
$\frac{F_{2}}{\mathrm{Im}\varphi^{F}}$ are Gorenstein flat, $\varphi^{F}$
is injective, we want to argue that $F$ is Gorenstein flat. By Lemma 3.8, both $\binom{F_{1}}{U\otimes F_{1}}$ and
$\binom{0}{\frac{F_{2}}{\mathrm{Im}\varphi^{F}}}$ are Gorenstein flat, and
there exists an short exact sequence in $T$-Mod
$$~~~0 \longrightarrow\binom{F_{1}}{U\otimes F_{1}}\longrightarrow\binom{F_{1}}{F_{2}}\longrightarrow\binom{0}{\frac{F_{2}}{Im\varphi^{F}}}\longrightarrow 0\\.$$
$T$ is coherent since $T$ is a Gorenstein ring, we know, by [17, Theorem 3.7], that $F=\binom{F_{1}}{F_{2}}\in\mathcal{GF}(T)$.\hfill$\Box$\\

In the following we want to characterize the
Gorenstein flat dimension of module over triangular matrix rings.
Using a similar way as we do in the proof of Theorem 2.5, we have the following result.
\begin{theorem} Let $n$ be a non-negative integer, $U_{A}$ and $_{B}U$ be finitely generated and have finite projective dimension, $T$ be a Gorenstein ring, $M=\binom{M_{1}}{M_{2}}_{\varphi^{M}}$ be a left $T$-module. Then $\mathrm{Gfd}(M)\leq n$ if and only if $\mathrm{Gfd}(M_{1})\leq n$,~$\mathrm{Gfd}(\frac{M_{2}}{\mathrm{Im}\varphi^{M}})\leq n$, and if $\binom{K_{1}}{K_{2}}_{\varphi^{K}}$ is a $n$-th syzygy of $M$, then $\varphi^{K}$ is injective.
\end{theorem}
\begin{proposition} Let $U_{A}$ and $_{B}U$ be finitely generated and have finite projective dimension, $T$ be Gorenstein. Let $M=\binom{M_{1}}{M_{2}}_{\varphi^{M}}$ be a left $T$-module. Then $M$, $M_{1}$ and $M_{2}$ have finite Gorenstein flat dimension.
\end{proposition}
\noindent{\bf{PROOF.}} Note that if $T$ is Gorenstein, every Gorenstein projective module is Gorenstein flat, so we have $\mathrm{Gfd}(M)\leq \mathrm{Gpd}(M)$. Since $T$ is Gorenstein, $T$ is left Gorenstein regular and $\mathrm{glGpd}(T)< \infty$. Thus $\mathrm{Gfd}(M)<\infty$.
By Lemma 3.2, $A$ and $B$ are Gorenstein, reasoning as before, $M_{1}$ and $M_{2}$ have finite Gorenstein flat dimension.\hfill$\Box$
\section{ Strongly Gorenstein homological dimensions}
~In this section, we present some characterizations of strongly Gorenstein projective (resp. injective, flat) modules over triangular matrix rings.

Recall that a left $R$-module $M$ is $strongly~Gorenstein~projective$ if there
exists an exact sequence of projective left $R$-modules
$\cdots \rightarrow P \stackrel{f}\rightarrow P \stackrel{f}\rightarrow  P \rightarrow \cdots$
with $M\cong\mathrm {Ker}f$ such that $\mathrm {Hom}_{R}(-,Q)$
leaves the sequence exact for any projective left
$R$-module $Q$. The strongly Gorenstein injective module are defined
dually. A left $R$-module $M$ is said to be $strongly~Gorenstein~flat$ if there exists an exact sequence $ \cdots \rightarrow F \stackrel{g}\rightarrow  F \stackrel{g}\rightarrow  F \rightarrow  \cdots$
of flat left $R$-modules with $M\cong\mathrm{Ker}g$ such that $I\otimes_{R}-$ leaves the sequence exact for any injective right $R$-module $I$. As usual, $\mathrm{SGpd}(M)$ and $\mathrm{SGid}(M)$ denote the strongly Gorenstein projective and injective dimensions of a left $R$-module $M$, respectively. Since every strongly Gorenstein projective left $R$-module is Gorenstein projective, by the proof of [5, Proposition 3.4] and [5, Theorem 3.5], we have the following conclusions.
\begin{proposition}Suppose that $U_{A}$ has finite flat dimension.

(1) If $M_{1}$ is a strongly Gorenstein projective left $A$-module, then $\mathbf{p}(M_{1},0)$ is a strongly Gorenstein projective left $T$-module.

(2) If $B$ is left Gorenstein regular, $_{B}U$ has
finite projective dimension and $(M_{1},M_{2})\in
A$-$\mathrm{Mod}\times B$-$\mathrm{Mod}$ is a strongly Gorenstein
projective object, then $\mathbf{p}(M_{1},M_{2})$ is a strongly
Gorenstein projective left $T$-module.
\end{proposition}
\noindent{\bf{PROOF.}} (1) Suppose that $M_{1}$ is strongly Gorenstein projective and let $\mathbf{P}:\cdots \rightarrow P \stackrel{f}\rightarrow P \stackrel{f}\rightarrow  P \rightarrow \cdots$ be an exact sequence consisting of projective left $A$-modules,
 which is $\mathrm{Hom}_{B}(-,C)$-exact for each projective left $A$-module $C$ and such that $\mathrm{Ker}\partial^{0}\cong M_{1}$. By [5, Lemma 2.3], we get that the complex $U\otimes_{A}\mathbf{P}$ is exact in $A$-Mod, which implies that the complex $\mathbf{p}(\mathbf{P})$ is exact in $T$-Mod. And it clearly verifies that $\mathrm{Ker}\partial^{0}_{\mathbf{p}(\mathbf{P})}=\mathbf{p}(M_{1},0)$. Finally, if $P=\binom{P_{1}}{P_{2}}$ is a projective left $T$-module, then the complex $\mathrm{Hom}_{T}(\mathbf{p}(\mathbf{P}),P)$ is isomorphic, by adjointness, to the
  complex $\mathrm{Hom}_{A}(\mathbf{P},P_{1})$, which is exact. This means that $\mathrm{Hom}_{T}(\mathbf{p}(\mathbf{P}),P)$ is exact, and so $\mathbf{p}(M_{1})$ is strongly Gorenstein projective.

(2) We only need to prove that both modules $\mathbf{p}(M_{1},0)$ and $\mathbf{p}(0,M_{2})$ are strongly Gorenstein projective when $M_{1}$ and $M_{2}$ are strongly Gorenstein projective. By (1), $\mathbf{p}(M_{1},0)$ is strongly Gorenstein projective.

Assume that $M_{2}$ is strongly Gorenstein projective and let $\mathbf{P}:\cdots \rightarrow P \stackrel{f}\rightarrow P \stackrel{f}\rightarrow  P \rightarrow \cdots$ be an exact sequence consisting of projective left $B$-modules which is $\mathrm{Hom}_{B}(_,C)$-exact for each projective left $B$-module and such that $\mathrm{Ker}\partial^{0}\cong M_{2}$. Then $\mathbf{p}(\mathbf{P})$ is an exact sequence of left $T$-modules such that $\mathrm{Ker}\partial^{0}\cong \mathbf{p}(0,M_{2})$. It remains to see that it is $\mathrm{Hom}_{T}(_,C)$-exact for each projective left $T$-module $C$. Let $C$ be a projective left $T$-module, and note that, as a consequence of [5, Corollary 2.3], there exists a projective object $(C_{1},C_{2})$ in $A$-Mod$\times$$B$-Mod such that $\mathbf{p}(C_{1},C_{2})=C$. Then, $C=\binom{C_{1}}{(U\otimes C_{1})\oplus C_{2}}$. Now, using adjointness, we get that the complex $\mathrm{Hom}_{T}(\mathbf{p}(\mathbf{P}),C)$ is isomorphic to the complex $\mathrm{Hom}_{B}(\mathbf{P},U\otimes C_{1})\oplus\mathrm{Hom}_{B}(\mathbf{P},C_{2})$. But $\mathrm{Hom}_{B}(\mathbf{P},C_{2})$ is exact, since $C_{2}$ is projective. In order to see that $\mathrm{Hom}_{B}(\mathbf{P},U\otimes_{A} C_{1})$ is exact, note that $U\otimes_{A} C_{1}$ has finite projective dimension in $B$-Mod, since it is isomorphic to a direct sum of copies of U, and using the condition that $B$ is left Gorenstein regular, it has finite injective dimension. Hence the exactness follows from [5, Lemma 2.4]. Consequently, $\mathrm{Hom}_{T}(\mathbf{P},C)$ is exact and the result is proved.\hfill$\Box$
\begin{theorem}
Let $U_{A}$ have finite flat dimension, $_{B}U$ have finite projective dimension and $B$ be left Gorenstein regular. Let $M=\binom{M_{1}}{M_{2}}_{\varphi^{M}}$  be a left $T$-module. Then the following are equivalent:

(1)~$M$ is strongly Gorenstein projective;

(2)~$M_{1}$ and $\frac{M_{2}}{\mathrm{Im}\varphi^{M}}$ are strongly Gorenstein projective and the $\varphi^{M}$ is injective.
\end{theorem}
\noindent{\bf{PROOF.}} (1)$\Rightarrow (2)$ Suppose $M$ is strongly Gorenstein projective in $T$-Mod
and let
$$~\mathbf{P}:~\cdots \longrightarrow
\binom{P_{1}}{P_{2}}\stackrel{\binom{\partial_{1}}{\partial_{2}}}
\longrightarrow\binom{P_{1}}{P_{2}}\stackrel{\binom{\partial_{1}}{\partial_{2}}}\longrightarrow
\binom{P_{1}}{P_{2}}\stackrel{}\longrightarrow
\cdots$$
be an exact sequence consisting of projective left $T$-modules
which is $\mathrm{Hom}_{T}(-,Q)$-exact for each projective left $T$-module $Q$ and
such that $\mathrm{Ker}\binom{\partial_{1}}{\partial_{2}}\cong M$. Then
$$\mathbf{P}_{1}:~~~~ \cdots \longrightarrow
P_{1}\longrightarrow P_{1}\stackrel{}\longrightarrow
P_{1}\longrightarrow\cdots~~~~~~~~~~\\$$
is an exact sequence consisting of projective left $A$-modules with $\mathrm{Ker}\partial_{1}\cong M_{1}$. Moreover, it is easy to see that it is $\mathrm{Hom}_{A}(-,Q)$-exact for each projective left $A$-module $Q$. This means that $M_{1}$ is a strongly Gorenstein projective left $A$-module.\\
\hspace*{6mm}Now note that for every morphism
$\binom{\partial_{1}}{\partial_{2}}:\binom{P_{1}}{P_{2}}\stackrel{}\longrightarrow
\binom{P_{1}}{P_{2}}$, we can construct the following
commutative diagram:
$$\CD
  0 @> >>U\otimes P_{1} @>\varphi^{P}>>P_{2}@>\pi>> \frac{P_{2}}{\mathrm{Im}\varphi^{P}}@> >> 0 \\
  @. @V 1_{U}\otimes\partial^{P}_{1} VV @V  \partial^{P}_{2} VV @V \overline{\partial^{P}_{2}} VV @.   \\
  0 @> >>U\otimes P_{1} @>\varphi^{P}>>P_{2}@>\pi>> \frac{P_{2}}{\mathrm{Im}\varphi^{P}}@> >> 0 ,
\endCD $$
where $\pi$ is the canonical projection. Using
this fact, the sequence $\mathbf{P}$ induce the sequence
$$\overline{\mathbf{P}_{2}}:~~~~ \cdots \longrightarrow
\frac{P_{2}}{\mathrm{Im}\varphi^{P}}\stackrel{\overline{\partial_{2}}}\longrightarrow
\frac{P_{2}}{\mathrm{Im}\varphi^{P}}\stackrel{\overline{\partial_{2}}}\longrightarrow
\frac{P_{2}}{\mathrm{Im}\varphi^{P}}\longrightarrow\cdots,~~~~~~~~~~\\$$
where $\varphi^{P}$ is the structural map of the $T$-module
$\binom{P_{1}}{P_{2}}$. Now each $\frac{P_{2}}{\mathrm{Im}\varphi^{P}}$ is projective in
$B$-Mod, since $\binom{P_{1}}{P_{2}}$ is a projective
left $T$-module. Moreover, the sequence
$\overline{\mathbf{P}_{2}}$ is exact, since there exists an short
exact sequence of complexes
$$~~~~0 \longrightarrow
U\otimes\mathbf{P}_{1}\longrightarrow
\mathbf{P}_{2}\stackrel{}\longrightarrow
\overline{\mathbf{P}_{2}}\longrightarrow0,~~~~~\\$$
in which $U\otimes\mathbf{P}_{1}$ is exact by [5, Lemma 2.3], and $\mathbf{P}_{2}$ is exact. It is easy to see $\mathrm{Ker}\overline{\partial_{2}}\cong\frac{M_{2}}{\mathrm{Im}\varphi^{M}}$. Since $B$ is left Gorenstein regular, so $\mathrm{Hom}_{B}(\overline{\mathbf{P}_{2}},Q)$ is exact for each projective left $B$-module $Q$ by [5, Lemma 2.4], which means $\frac{M_{2}}{\mathrm{Im}\varphi^{M}}$ is strongly Gorenstein projective.

Finally, we prove that morphism $\varphi^{M}$ is
injective. By [5, Lemma 2.3], $U\otimes\mathbf{P}_{1}$ is an exact
sequence. This means that if $\iota_{1}:M_{1}\rightarrow P_{1}$
is the inclusion, $1_{U}\otimes\iota_{1}$ is monic. But, since
$\binom{M_{1}}{M_{2}}_{\varphi^{M}}$ is a submodule of
$\binom{P_{1}}{P_{2}}_{\varphi^{P}}$, the following diagram
commutes:
$$\CD
 U\otimes M_{1} @>1_{U}\otimes\iota_{1}>>U\otimes P_{1}  \\
@V\varphi^{M} VV @V \varphi^{P} VV   \\
  M_{2} @>\iota_{2}>> P_{1},
\endCD $$
where $\iota_{2}$ is the inclusion. Since $\varphi^{P}$ is monic, and $\binom{P_{1}}{P_{2}}_{\varphi^{P}}$ is strongly Gorenstein projective, then we conclude that $\varphi^{M}$ is injective.

(2)$\Rightarrow$ (1) Assume $M_{1}$ and $\frac{M_{2}}{\mathrm{Im}\varphi^{M}}$ are strongly Gorenstein projective and $\varphi^{M}$ is injective.
Then we have a exact sequence $\mathbf{Q}:\cdots\rightarrow Q\stackrel{d'}\rightarrow Q\stackrel{d'}\rightarrow Q\rightarrow \cdots$ of projective $A$-modules, which is $\mathrm {Hom}_{A}(-,C)$-exact for any projective left $A$-module $C$ and such that $M_{1}\cong\mathrm{Ker}d'$. By [5, Lemma 2.3], $U\otimes \mathbf{Q}$ is exact, and so $0\rightarrow U\otimes M_{1}\stackrel{}\rightarrow U\otimes Q\stackrel{1_{U}\otimes d'}\rightarrow U\otimes Q\rightarrow \cdots$ is exact. Since $\frac{M_{2}}{\mathrm{Im}\varphi^{M}}$ is strongly Gorenstein projective, we have a exact sequence $\mathbf{P}:\cdots\rightarrow P\stackrel{d}\rightarrow P\stackrel{d}\rightarrow P\rightarrow \cdots$ of projective $B$-modules, which is $\mathrm {Hom}_{B}(-,C)$-exact for any projective left $B$-module $C$ and such that $\frac{M_{2}}{\mathrm{Im}\varphi^{M}}\cong\mathrm{Ker}d$. Since $_{B}U$ has finite projective dimension, it follows from [3, Proposition 2.9] that $\mathrm{Ext}^{1}_{B}(\mathrm{Ker}d,U)=0$. Since $Q$ is a projective $A$-module, $\mathrm{Ext}^{1}_{B}(\mathrm{Ker}d,U\otimes Q)=0$. Apply [19, Lemma 1.6(1)] to exact sequence
$0 \longrightarrow U\otimes M_{1}\longrightarrow M_{2}\stackrel{}\longrightarrow \frac{M_{2}}{\mathrm{Im}\varphi^{M}}\longrightarrow0$, we obtain an exact sequence
$$0\longrightarrow M_{2}\longrightarrow (U\otimes Q)\oplus P\stackrel{\partial}\longrightarrow (U\otimes Q)\oplus P\longrightarrow \cdots\eqno(1)$$
with $\partial=\left(\begin{smallmatrix}  d & 0 \\  \sigma & 1\otimes d' \\  \end{smallmatrix}\right)$, $\sigma:P\rightarrow U\otimes M$. Apply Generalized Horseshoe Lemma [19, Lemma 1.6] to exact sequence $0 \longrightarrow U\otimes M_{1}\longrightarrow M_{2}\stackrel{}\longrightarrow \frac{M_{2}}{\mathrm{Im}\varphi^{M}}\longrightarrow0$, the exact sequence $\cdots\rightarrow U\otimes Q\stackrel{1_{U}\otimes d'}\rightarrow U\otimes Q\rightarrow U\otimes M_{1}\rightarrow0$ and $\cdots\rightarrow P\stackrel{d}\rightarrow P\stackrel{}\rightarrow \frac{M_{2}}{\mathrm{Im}\varphi^{M}}\rightarrow 0$, we obtain another exact sequence
$$\cdots \longrightarrow (U\otimes Q)\oplus P\stackrel{\partial}\longrightarrow (U\otimes Q)\oplus P\longrightarrow M_{2}\longrightarrow 0.\eqno(2)$$
Putting (1) and (2) together we get the following exact sequence of projective $T$-modules
$$\mathbf{L}:\cdots\longrightarrow \binom{Q}{(U\otimes Q)\oplus P}\stackrel{\binom{d'}{\partial}}\longrightarrow\binom{Q}{(U\otimes Q)\oplus P}\stackrel{\binom{d'}{\partial}}\longrightarrow\binom{Q}{(U\otimes Q)\oplus P}\longrightarrow \cdots$$
with $\mathrm{Ker}\binom{d'}{\partial}\cong M$. In fact, $\mathbf{L}=\mathbf{p}(\mathbf{Q},\mathbf{P})$. Let $C$ be a projective left $T$-module. As a consequence of [5, Corollary 2.3], there exists a projective object $(C_{1},C_{2})$ in $A$-Mod$\times$$B$-Mod such that $\mathbf{p}(C_{1},C_{2})=C$. Then, $C=\binom{C_{1}}{(U\otimes C_{1})\oplus C_{2}}$. Now, using adjointness, we get that the complex $\mathrm{Hom}_{T}(\mathbf{p}(\mathbf{Q},\mathbf{P}),C)$ is isomorphic to the complex $\mathrm{Hom}_{B}(\mathbf{Q},C_{1})\oplus\mathrm{Hom}_{B}(\mathbf{P},U\otimes C_{1})\oplus\mathrm{Hom}_{B}(\mathbf{P},C_{2})$. But $\mathrm{Hom}_{B}(\mathbf{Q},C_{1})$ and $\mathrm{Hom}_{B}(\mathbf{P},C_{2})$ are exact, since $C_{1}$ and $C_{2}$ is projective. In order to see that $\mathrm{Hom}_{B}(\mathbf{P},U\otimes_{A} C_{1})$ is exact, note that $U\otimes_{A} C_{1}$ has finite projective dimension in $B$-Mod, since it is isomorphic to a direct sum of copies of $U$, and using the condition that $B$ is left Gorenstein regular, it has finite injective dimension. Then the exactness follows from [5, Lemma 2.4]. Consequently, $\mathrm{Hom}_{T}(\mathbf{L},C)$ is exact and the result is proved.\hfill$\Box$\\

By analogous arguments, we have the following results.
\begin{proposition}Suppose that $_{B}U$ has finite flat dimension.

(1) If $M_{2}$ is a strongly Gorenstein injective left $B$-module, then $\mathbf{h}(0,M_{2})$ is a strongly Gorenstein injective left $T$-module.

(2) If $A$ is left Gorenstein regular, $U_{A}$ has
finite flat dimension and $(M_{1},M_{2})\in A$-$\mathrm{Mod}\times
B$-Mod is a strongly Gorenstein injective object, then
$\mathbf{h}(M_{1},M_{2})$ is a strongly Gorenstein projective left
$T$-module.
\end{proposition}
\begin{proposition}
Suppose that $U_{A}$ has finite flat dimension, $_{B}U$ has finite projective dimension and $A$ is left Gorenstein regular. Let $M=\binom{M_{1}}{M_{2}}_{\varphi^{M}}$  be a left T-module. Then the following are equivalent:

$(1)~M$ is strongly  Gorenstein injective;

$(2)~M_{2}$ and $\mathrm{Ker}\widetilde{\varphi^{M}}$ are strongly  Gorenstein injective and the $\widetilde{\varphi^{M}}$ is a surjection.
\end{proposition}
\begin{proposition} Let $U_{A}$ and $_{B}U$ have finite projective dimension, $T$ be a right Gorenstein regular ring, $M=\binom{M_{1}}{M_{2}}_{\varphi^{M}}$ a left $T$-module. If $M$ is strongly Gorenstein flat, then $M_{1}$ and $\frac{M_{2}}{\mathrm{Im}\varphi^{M}}$ are strongly Gorenstein flat and the morphism $\varphi^{M}$ is injective.
\end{proposition}
\noindent{\bf{PROOF.}} It is similar to the proof of Proposition 3.6.
\hfill$\Box$\\

One can prove the following argument in a similar way as we do in the proof of Theorem 2.5.
\begin{theorem} Let $n$ be a non-negative integer, $U_{A}$ have finite flat dimension, $_{B}U$ have
finite projective dimension and $B$ be left Gorenstein regular. Let
$M=\binom{M_{1}}{M_{2}}_{\varphi^{M}}$  be a left $T$-module. Then
$\mathrm{SGpd}(M)\leq n$ if and only if $\mathrm{SGpd}(M_{1})\leq n$,
$\mathrm{SGpd}(\frac{M_{2}}{\mathrm{Im}\varphi^{M}})\leq n$ and if
$\binom{K_{1}}{K_{2}}_{\varphi^{K}}$ is a n-th syzygy of $M$, then
$\varphi^{K}$ is injective.
\end{theorem}
\begin{theorem} Let $n$ be a non-negative integer, $U_{A}$ have finite flat dimension, $_{B}U$ have finite projective dimension and $A$ be left Gorenstein regular. Let $M=\binom{M_{1}}{M_{2}}_{\varphi^{M}}$  be a left $T$-module. Then $\mathrm{SGid}(M)\leq n$ if and only if  $\mathrm{SGid}(M_{2})\leq n$,~$\mathrm{SGid}\mathrm{(Ker}\widetilde{\varphi^{M}})\leq n$ and if $\binom{L_{1}}{L_{2}}_{\varphi^{L}}$ is a n-th cosyzygy of $M$, then $\widetilde{\varphi^{L}}$ is surjective.
\end{theorem}

\end{document}